\documentclass[12pt,letterpaper,reqno]{amsart}
\usepackage[utf8]{inputenc}
\usepackage[T1]{fontenc}
\usepackage{amsmath}
\usepackage{amsfonts}
\usepackage{amssymb}
\usepackage{graphicx}
\usepackage{xcolor}
\usepackage{subfig}

\usepackage{amsmath,amssymb,latexsym,graphicx}

\allowdisplaybreaks

\newtheorem{theorem}{Theorem}

\allowdisplaybreaks

\begin{document}
	
\author{Sergii Myroshnychenko, Kateryna Tatarko, and Vladyslav Yaskin}

\address{S.~Myroshnychenko, Department of Mathematical Sciences, Lakehead University, Barrie, ON L4M 3X9, Canada}

\email{smyroshn@lakeheadu.ca}

\address{K.~Tatarko, Department of Pure Mathematics, University of Waterloo, Waterloo, ON N2L 3G1, Canada}

\email{ktatarko@uwaterloo.ca}

\address{V.~Yaskin, Department of Mathematical and Statistical Sciences, University of Alberta, Edmonton, AB T6G 2G1, Canada}

\email{yaskin@ualberta.ca}

 \thanks{The second and third  authors were supported in part by NSERC. The second and third authors are  also grateful to ICERM for excellent working conditions, where they participated  in the semester-long program ``Harmonic Analysis and Convexity''.}

 
	
\title[How far apart can centroids be?]{How far apart can the projection of the centroid of a convex body and the centroid of its projection be?}

\begin{abstract} We show that there is a constant $D \approx 0.2016$ such that for every $n$, every convex body $K\subset \mathbb R^n$, and every hyperplane $H\subset \mathbb R^n$, the distance between the projection of the centroid of $K$ onto $H$ and the centroid of the projection of $K$ onto $H$ is at most $D$ times the width of $K$ in the direction of the segment connecting the two points. The constant $D$ is asymptotically sharp.
	\end{abstract}

\maketitle

\section{Introduction}

Let $K$ be a convex body in $\mathbb R^n$, i.e., a compact convex set with non-empty interior. The centroid (the center of mass) of $K$ is the point 
$$c(K) = \frac{1}{|K|} \int_K x \, dx,$$
where $|K|$ denotes the volume of $K$ and the integration is with respect to Lebesgue measure.

In this paper we study the following  question. Let $H$ be a hyperplane in $\mathbb R^n$. Denote by $P_Hc(K)$  the orthogonal projection of the centroid of $K$ onto $H$ and by $c(P_HK)$  the centroid of the projection of $K$ onto $H$. For centrally symmetric bodies these two points coincide, but for non-symmetric bodies these points are generally different. Thus it is natural to ask how far apart these two points can be relative to some linear size of $K$. More precisely, we are interested in the smallest constant $D_n$ such that for any convex body $K$ in $\mathbb R^n$ we have $$ |P_Hc(K)-c(P_HK)|\le D_n w_K (u),$$ where $u$ is the unit vector parallel to the segment connecting $ P_Hc(K)$ and $c(P_HK)$, and $w_K(u)$ is the width of $K$ in the direction of $u$, given by 
$$w_K(u) = \max_{x\in K} \{\langle x,u\rangle\} -  \min_{x\in K} \{\langle x,u\rangle\}.$$
Questions of this type began attracting attention several years ago in connection to    Gr\"unbaum-type inequalities for sections and projections;    see \cite{SZ}, \cite{MSZ}. In particular, an analogue of the question  above for sections of convex bodies is stated in  \cite[p.~127]{S}. 

For other questions related to distances between various centroids the reader is referred to the book \cite[p.~36]{CFG} and the references contained therein.

\section{Main results}\label{main}
Let us start with the simplest case when $n=2$. The projection of $K$ on any line is a segment whose centroid is just its midpoint. On the other hand, if we have two parallel supporting lines to $K$ then the distance from the centroid of $K$ to one these lines is at least $1/3$ of the distance between the lines; see \cite[p.~58]{BF}. Thus  in the 2-dimensional case we have
 $$ |P_Hc(K)-c(P_HK)|\le \frac16 w_K (u).$$
The constant $1/6$ is sharp; it is attained on triangles with one side perpendicular to the line $H$. 

We will now present higher-dimensional analogues of this observation. 

\begin{theorem} Let $D_n$, $n\ge 3$,  be the smallest number such that 
	\begin{equation}\label{MainIneq}
	|P_Hc(K)-c(P_HK) |\le D_n\cdot w_K (u),
	\end{equation}
for every convex body $K$  in $\mathbb R^n$ and   every hyperplane $H\subset \mathbb R^n$,
	where $$u = \frac{P_Hc(K)-c(P_HK)}{|P_Hc(K)-c(P_HK)|},$$
	provided $P_Hc(K)\ne c(P_HK)$. Then

\begin{enumerate}\item[(i)]	  $D_3 = 1-\displaystyle\sqrt{\frac23}\approx 0.1835$;  	 the sequence  $\{D_n\}_{n=3}^\infty$ is increasing;  and  $\displaystyle\lim_{n\to \infty} D_n\approx  0.2016.$
	
\item[(ii)]		Inequality (\ref{MainIneq}) turns into equality if and only if $K$ is a body obtained as follows. For a fixed hyperplane  $H$ and a vector $u$ parallel to $H$, denote by $\theta$ a unit normal vector to $H$ and take any $(n-2)$-dimensional subspace $U$ orthogonal to $u$ and transversal to $\theta$. Let 
  $L_0$ be any convex body in $U$. Denote by $tL_0$  the dilation of $L_0$ with respect to its centroid by a factor of $t=t_{max}$, which will be defined later in the proof.
Let $\lambda$, $\mu$, $\nu$ be real numbers, $\mu\ne \nu$.   Define $L_1= tL_0+\lambda u+\mu \theta$ and $L_2= tL_0+\lambda u+\nu \theta$. Then $K$ is the convex hull of $L_0$, $L_1$, and $L_2$. Figure \ref{fig:optimal} shows an example of such a body in $\mathbb R^3$ when $H=\{x_3=0\}$, $u=e_1$, and $U$ is the linear span of $e_2$.

\end{enumerate}
  \begin{figure}[h]
  	\includegraphics[scale=0.5]{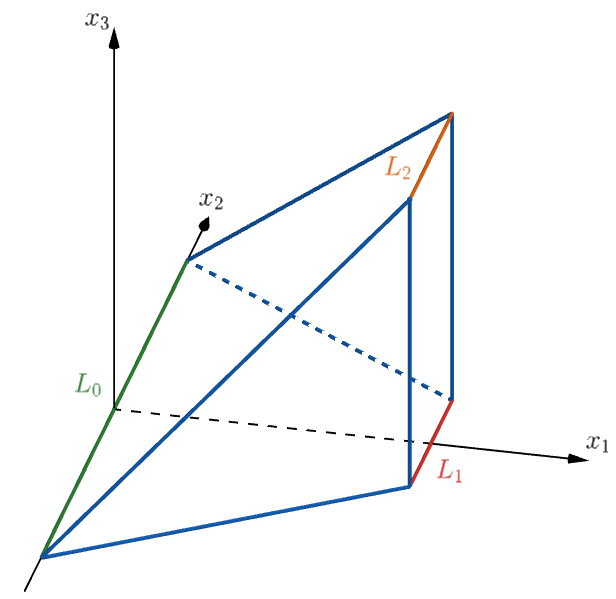}
  	\caption{
  	}
  	\label{fig:optimal}
  \end{figure}

	\end{theorem}

  \begin{proof} The proof will consist of several steps. 
  	
  	{\it Step 1 - Reducing the problem to a smaller class of bodies.}  	
  	 Without loss of generality we may assume that the hyperplane $H$ is the coordinate plane $x_n=0$. Moreover, it is enough to prove the theorem in the class of convex bodies $K$ described below. First, we will assume that the points $P_Hc(K)$ and $c(P_HK)$ lie on the $x_1$-axis, and the supporting hyperplanes to $K$ perpendicular to $e_1$ are the planes $x_1=0$ and $x_1=1$. Note that for such bodies $w_K(e_1)=1$.
  	 Additionally, we may assume that $\langle P_Hc(K),e_1\rangle\ge \langle c(P_HK),e_1\rangle.$ 
  	
  	 Next let us  perform the following operation that does not change $P_Hc(K)$ and $P_H(K)$. For each line $\ell$ orthogonal to $x_n=0$, replace the segment $K\cap \ell$  with another segment of the same length on the line $\ell$ that lies in the half-space $x_n\ge 0$, with one end-point on the plane $x_n=0$. This process is sometimes called Blaschke shaking; see \cite[p. 92]{Ga} for details. 
  	   	 	Now we will apply the following Schwarz-Steiner symmetrization. We transform $K$ into a convex body $\widetilde K$ such that for each $(n-2)$-dimensional affine subspace $U$ orthogonal to the two-dimensional subspace $V=\mbox{span}\{e_1,e_n\}$ the section $\widetilde K\cap U$ is an $(n-2)$-dimensional Euclidean ball  centered at the point $U\cap V$ with the same $(n-2)$-dimensional volume as $ K\cap U$. Abusing notation we will still denote the new body by $K$. Observe that the operation performed above does not change $P_Hc(K)$ and $c(P_HK)$.   	 
  	 Also note that $P_H K$ is a body of revolution in $H$ with $x_1$ as its axis of revolution.
  	 
  	 Thus we have reduced the problem to the class of convex bodies of the form
   $K=\{x\in \mathbb R^n:  0\le x_1\le 1,\, \sqrt{x_2^2+\cdots+x_{n-1}^2}\le f(x_1), 0\le x_n\le g(x_1,\sqrt{x_2^2+\cdots+x_{n-1}^2})\}$, where  $f$ and $g$ are concave functions.   	 
   Without loss of generality, we can rescale the functions $f$ and $g$ so that their maximal values are equal to 1. 
   The set of convex bodies introduced above will be denoted by $\mathcal C$. Then 
\begin{equation}\label{sup-in-C}
D_n=\sup_{K\in \mathcal C} \langle  c(K)-c(P_HK),e_1\rangle .
\end{equation}

   Observe that all bodies from $\mathcal C$   
     are contained in the box $\{0\le x_1\le 1,|x_2|\le 1, \ldots, |x_{n-1}|\le 1, 0\le x_n\le 1\}$.   Moreover, their volumes are bounded from below by $\frac1{n(n-1)} \kappa_{n-2},$ where $\kappa_{n-2}$ is the volume of the $(n-2)$-dimensional Euclidean ball. This is easily seen by observing that the height of the bodies is 1 (in the direction $x_n$)  and  $|P_HK|\ge \frac1{n-1} \kappa_{n-2}$.
      Here and below we use $|\cdot|$ to denote the $k$-dimensional  volumes of  bodies in $k$-dimensional subspaces.

By the Blaschke selection theorem \cite[p.~63]{Sch} there is a convex body for which the supremum   
 in (\ref{sup-in-C}) is achieved. 
  	 Moreover, such a body has non-empty interior, since the volumes of the bodies in the class  $\mathcal C$ are bounded below by a positive constant. 
  	
  	Throughout the text we will call maximizers those bodies for which the supremum in (\ref{sup-in-C}) is achieved (or more generally, those convex bodies for which there is equality in (\ref{MainIneq})).

 {\it Step 2 - Describing the projection of maximizers.} Let $K\in \mathcal C$ be a maximizer.	
   We claim  that $P_H K$ is obtained by rotating a trapezoid in $H$ around the $x_1$-axis.    
To this end, let us show that there is $\xi \in {S}^1  \subset \text{span}\{e_1, e_2\}$ and an affine hyperplane $\beta$ in $\mathbb R^n$ orthogonal to $\xi$ 
   such that $\beta$ divides both $K$ and $P_HK$ into two subsets of equal volumes. To see this, first observe that for every $\xi \in {S}^1 =S^{n-1}\cap   \text{span}\{e_1, e_2\}$ there is an
     affine hyperplane $\beta_{\xi}$ orthogonal to $\xi$  that divides $P_HK$ into two parts of equal $(n-1)$-dimensional volumes. Suppose $\beta_\xi$ is given by the equation $\langle x,\xi\rangle = s(\xi)$, for some function $s$.     Consider
  the function
  $$
  \Phi(\xi) = |K \cap \beta^+_{\xi}| - |K \cap \beta^-_{\xi}|, \quad \xi\in S^1,
  $$
  where $\beta^+_{\xi} = \{x\in\mathbb R^n: \langle x,\xi\rangle \ge s(\xi)\}$ and $\beta^-_{\xi} = \{x\in\mathbb R^n: \langle x,\xi\rangle \le s(\xi)\}$.
  Since $\Phi$ is continuous and odd, there is $\xi_0\in S^1$, such that $\Phi(\xi_0) = 0$. Thus $\beta=\beta_{\xi_0}$ is the required affine hyperplane.  Denote by $K_1^1$ and $K_2^1$ the two subsets of equal volume into which $\beta$ splits $K$. Then
  $$
  c(K) = \frac{1}{|K|} \left(c(K_1^1) |K_1^1| + c(K_2^1)|K_2^1|\right) = \frac{1}{2}  \left(c(K_1^1)   + c(K_2^1) \right),
  $$
  and similarly,
    $$
  c(P_HK) = \frac{1}{|P_HK|} \left(c(P_HK_1^1) |P_HK_1^1| + c(P_HK_2^1)|P_HK_2^1|\right) = \frac{1}{2}  \left(c(P_HK_1^1)   + c(P_HK_2^1) \right).
  $$
Therefore,
  \begin{align*}
  D_n &= \langle c(K) - c(P_HK), e_1\rangle \\ &= \frac{1 }{2}\langle c(K_1^1)-c(P_HK_1^1)+c(K_2^1) - c(P_HK_2^1), e_1\rangle  \\
  &\le \frac{1}{2} \left(D_n \cdot w_{K_1^1}(e_1)+D_n \cdot w_{K_2^1}(e_1)\right)\le D_n \cdot w_K(e_1)=D_n.
  \end{align*}
This means that the above inequalities are in fact equalities, and  
  $$
  w_{K_1^1}(e_1) = w_{K_2^1}(e_1)= w_K(e_1) =1.
  $$
  Moreover,  
  $$
 \langle c(K_1^1)-c(P_HK_1^1), e_1\rangle  = \langle  c(K_2^1) - c(P_HK_2^1), e_1\rangle =D_n.
  $$
  In other words, $K_1^1$ and $K_2^1$ are also maximizers. Thus, we can apply the same procedure to  the bodies $K_1^1$ and $K_2^1$, and produce four maximizers $K_1^2$, $K_2^2$, $K_3^2$, and $K_4^2$, and so on.  After doing this $m$ times, we will get $2^m$ convex bodies $\{K_i^m\}_{i=1}^{2^m}$ of volume $2^{-m}|K|$, and the width of each such body in the direction $e_1$ equals 1. Note that all these bodies are maximizers. The set of planes used to obtain these $2^m$ convex bodies will be denoted by $\Pi_m$. Suppose that $P_HK$ is generated  by rotating the figure $\{(x_1,x_2): 0\le x_1\le 1, \, 0\le x_2\le f(x_1)\}$ around the $x_1$-axis in $H$. To reach a contradiction, assume that $f$ is not a linear function. Let $\langle x,\xi\rangle = \sigma$ be the line in the $x_1x_2$-plane that passes through the points $(0,f(0))$ and $(1,f(1))$, and let the curve $x_2=f(x_1)$ be contained in the half-plane  $\{\langle x,\xi\rangle \ge \sigma\}$. 
  This half-plane is the intersection of the $x_1x_2$-plane with  the half-space $\{x\in \mathbb R^n: \langle x,\xi\rangle \ge \sigma\}$.  If $f$ is not a linear function, then the set $A=K\cap \{x\in \mathbb R^n: \langle x,\xi\rangle \ge \sigma\}$ has strictly positive  volume.  On the other hand none of the planes from the set $\Pi_m$ can intersect the interior of $A$; otherwise the width of one of the   bodies $K_i^m$ would be strictly less than 1. Therefore for each $m$ there is a body $K_i^m$ that contains $A$, but this is impossible since the volume of $A$ is fixed and the volume of $K_i^m$ is $2^{-m}|K|$. We reached a contradiction, which means that the function $f$ generating $P_HK$ is linear.

 {\it Step 3 -  Describing the upper part of the boundary of maximizers.} 	
  	As was shown above, we can assume that the projection of a maximizer onto $H$ is of the form $\{(x_1,...,x_{n-1}):  0\le x_1\le 1,\, \sqrt{x_2^2+\cdots+x_{n-1}^2}\le ax_1+b\} $, where $a$ and $b$ are constants.
  	 Denote the latter convex set in $H$ by $M$ and consider two convex bodies in $\mathbb R^n$: $K=\{x\in \mathbb R^n: (x_1,...,x_{n-1})\in M, 0\le x_n\le g(x_1,\sqrt{x_2^2+\cdots+x_{n-1}^2})\}$ and $L=\{x\in \mathbb R^n: (x_1,...,x_{n-1})\in M, 0\le x_n\le x_1\}$. We will show that $\langle  c(K),e_1\rangle \le \langle  c(L),e_1\rangle$.

  	 Let us first compute the volumes of  $K$ and $L$.
  	 \begin{align*}
  	 |K| &=  \int_{M} g\left(x_1,\sqrt{x_2^2+\cdots+x_{n-1}^2}\right) \,  dx_1\ldots dx_{n-1} \\
  	 & = 
  	  \int_0^1 \int_{{S}^{n-3}} \int_0^{a   x_1 + b}  r^{n-3} g(x_1,r) dr\, d\theta \, dx_1\\ &= |\mathbb{S}^{n-3}| \int_0^1 \int_0^{a   x_1 + b} r^{n-3} g(x_1,r)\, dr\, dx_1 = \\
  	 & =|\mathbb{S}^{n-3}| \int_0^1(a   x_1 + b)^{n-2} \int_0^1 z^{n-3} g(x_1,(a  x_1 + b) z) \, dz \, dx_1  \\
  	 &=|\mathbb{S}^{n-3}| \int_0^1(a   x_1 + b)^{n-2} \cdot G(x_1) \, dx_1,
  	 \end{align*}
  	 where $$
  	 G(x_1) = \int_0^1 z^{n-3} g(x_1,(a   x_1+b)z) \, dz.
  	 $$
  	 
  	 Since $g(x_1,r )$ is a concave function of two variables, it follows that $G$ is   a concave function of one variable. Indeed, for any $\tau_1$, $\tau_2\in [0,1]$ we have
  	 \begin{align}\label{G_concave}
  	 G\left(\frac{\tau_1+\tau_2}{2}\right) =& \int_0^1 z^{n-3} g\left(\frac{\tau_1+\tau_2}{2}, \left(a \frac{\tau_1+\tau_2}{2}+b\right)z\right) \, dz\nonumber\\
  	 &=\int_0^1 z^{n-3} g\left(\frac{\tau_1+\tau_2}{2}, \frac{(a   \tau_1 +b)z+(a  \cdot \tau_2+b)z}{2}\right) \, dz \nonumber \\
  	 &\geq\int_0^1 z^{n-3} \cdot \frac{g(\tau_1,(a   \tau_1+b)z)+g(\tau_2,(a   \tau_2+b)z)}{2} \, dz  \\
  	 & =\frac{G(\tau_1)+G(\tau_2)}{2}. \nonumber
  	 \end{align}
  	 
  	Similarly, we have
  	 \begin{align*}
  	 |L|& =    \int_{M} x_1 \, dx_1\ldots \, dx_{n-1} =|\mathbb{S}^{n-3}| \int_0^1 \int_0^{a   x_1 + b} r^{n-3} \cdot x_1 \, dr\, dx_1\\
  	 & = |\mathbb{S}^{n-3}| \int_0^1 \frac{ (a \cdot x_1 + b)^{n-2} \cdot x_1}{n-2} \, dx_1.
  	 \end{align*}

Multiplying $g$ by a constant, if needed, we may assume that $|K| = |L|$, which is equivalent to
$$
\int_0^1 (a   x_1 + b)^{n-2} \left(G(x_1) - \frac{x_1}{n-2}\right) \, dx_1 = 0.
$$
The latter implies that there is $\alpha \in (0,1)$ such that
 $G(\alpha) = \frac{\alpha}{n-2}$. Moreover, by the concavity of $G$ such $\alpha$ is unique (unless $G$ is itself linear). Hence, 
 $G(x_1) \ge \frac{x_1}{n-2},$ when $ 0\le  x_1\le \alpha$, and
 $G(x_1) \le \frac{x_1}{n-2},$ when $ \alpha\le  x_1\le 1.$
 
Therefore,
\begin{align*}
&\int_0^1 (a   x_1 + b)^{n-2} \cdot \left(G(x_1) - \frac{x_1}{n-2}\right) \cdot x_1 \, dx_1 = \\
&\int_0^1 (a  x_1 + b)^{n-2} \cdot \left(G(x_1) - \frac{x_1}{n-2}\right) \cdot (x_1 - \alpha) \, dx_1 \leq 0,
\end{align*}
meaning that $$\int_K x_1\, dx \le \int_L x_1\, dx,$$
i.e., $\langle  c(K),e_1\rangle \le \langle  c(L),e_1\rangle$.
The latter inequality is strict, unless $G(x_1) = \frac{x_1}{n-2}$.
In Step 7 we will show that this implies  $g(x_1,r)=x_1$.

{\it Step 4 - Describing $a$ and $b$ for maximizers.}
As we have shown, we need to look for maximizers among the bodies of the form $$K=\{x\in\mathbb R^n:  0\le x_1\le 1,\, \sqrt{x_2^2+\cdots+x_{n-1}^2}\le ax_1+b,\, 0\le x_n\le x_1\}.$$ It remains to determine what constants   $a$ and $b$ yield the maximum. Since $ax_1+b\ge 0$ for $0\le x_1\le 1$, we necessarily have $b\ge 0$ and $a+b\ge 0$.

One can see that
\begin{align*}
\langle	c(K) - c(P_HK), e_1\rangle &=
	\frac{\int_0^1 (a  x_1 + b)^{n-2}   x_1^2 \, dx_1}{ \int_0^1(a   x_1 + b)^{n-2}   x_1 \, dx_1} -\frac{\int_0^1 (a   x_1 + b)^{n-2}   x_1 \, dx_1}{ \int_0^1(a   x_1 + b)^{n-2} \, dx_1}\\
	&= \frac{\int_0^1 (t   x_1 + 1)^{n-2}   x_1^2 \, dx_1}{ \int_0^1(t   x_1 + 1)^{n-2}   x_1 \, dx_1} -\frac{\int_0^1 (t   x_1 + 1)^{n-2}   x_1 \, dx_1}{ \int_0^1(t   x_1 + 1)^{n-2} \, dx_1}
\end{align*}
where $t = \frac ab\ge -1$, and $b=0$ corresponds to $t\to \infty$.

Computing the above integrals, we get
\begin{align*} 
\int_0^1 (t x_1 + 1)^{n-2} \, dx_1  &= \frac{(t+1)^{n-1}-1}{(n-1)t},\\
\int_0^1 (t x_1 + 1)^{n-2} \cdot x_1 \, dx_1 & = \frac{1+(t+1)^{n-1}((n-1)t-1)}{n(n-1)t^2},
\\
	\int_0^1 (t x_1 + 1)^{n-2} \cdot x_1^2 \, dx_1 & = \frac{-2+(t+1)^{n-1} (2+(n-1)t(nt-2))}{n(n^2-1)t^3},
\end{align*}
and hence,
\begin{align*} 
&\langle c(K) - c(P_HK),e_1\rangle\\
&\qquad=
\frac{-2+(t+1)^{n-1} (2+(n-1)t(nt-2))}{(n+1)t(1+(t+1)^{n-1}((n-1)t-1))} -  \frac{1+(t+1)^{n-1}((n-1)t-1)}{n t((t+1)^{n-1}-1)}.
\end{align*}
Let us denote the latter expression by $D_n(t)$. Thus, $D_n=\sup_{t\ge -1} D_n(t)$.

{\it Step 5 - Finding $D_3$.} Computing the supremum of $D_n(t)$ cannot be done precisely   in all dimensions, since this requires solving polynomial equations of high degrees. However, in the three-dimensional case this is doable. 
When $n = 3$, we get 
\begin{align*}
D_3(t)& = \frac{-2+(t+1)^{2} (2+2t(3t-2))}{4t(1+(t+1)^2(2t-1))} -  \frac{1+(t+1)^2(2t-1)}{3 t((t+1)^2-1)}\\
 &=\frac{3t+4}{2(2t+3)}-\frac{2t+3}{3(t+2)}=\frac{1}{6} \cdot \frac{t^2+6t+6}{2t^2+7t+6}.
\end{align*}
To find the maximum of $D_3(t) $, let us compute its derivative,
$$
D_3'(t) = -  \frac{1}{6} \cdot \frac{5t^2 + 12 t + 6}{(2t^2+7t+6)^2}.
$$
The only critical number of $D_3(t) $ in the interval  $[-1, \infty)$ is $t = \frac{\sqrt{6} - 6}5$. This corresponds to the point of global maximum of $D_3(t)$ on $[-1, \infty)$, since $D_3'(t)>0$ when $-1< t<\frac{\sqrt{6} - 6}5$ and $D_3'(t)<0$ when $t>\frac{\sqrt{6} - 6}5$. Thus,   $D_3 = 1-\sqrt{\frac23}\approx 0.1835$ which is attained at $t =   \frac{\sqrt{6} - 6}5$.

\medskip

{\it Step 6 - Finding $\sup_{n} D_n$.}
Let us now discuss the case of higher dimensions. First, we will show that 
the sequence $\{D_n\}_{n=3}^\infty$ is non-decreasing. Let $K\subset \mathbb R^n$ be a maximizer and  $K\in \mathcal C$. Consider a convex body $\widetilde K\subset \mathbb R^{n+1}$ defined by $$ \widetilde K= K \times [0,1] =\{(x_0,x_1,\ldots, x_n): 0\le x_0\le 1, (x_1,\ldots, x_n)\in K\}.$$
As before, let $H$ be the hyperplane in $\mathbb R^n$ defined by $\{(x_1,\ldots,x_n):x_n=0\}$, and $\widetilde H$ be the hyperplane in $\mathbb R^{n+1}$ defined by $\{(x_0,x_1,\ldots,x_n):x_n=0\}$.
Then
$$ D_n = \langle c(K) - c(P_HK),e_1\rangle = \langle c(\widetilde K) - c(P_{\widetilde H}\widetilde K),e_1\rangle \le D_{n+1} w_{\widetilde K}(e_1)=D_{n+1}, $$
as claimed.

We will now compute the limit $\displaystyle \lim_{n\to \infty} D_n = \sup_n D_n$.
Recall that 
\begin{align*}
D_n(t) = \frac{\int_0^1 (t x_1 + 1)^{n-2} x_1^2 dx_1}{\int_0^1 (t x_1 + 1)^{n-2} x_1 dx_1} - \frac{\int_0^1 (t x_1 + 1)^{n-2} x_1 dx_1}{\int_0^1 (t x_1 + 1)^{n-2} dx_1}.
\end{align*}
Our goal is to show that the latter quantity will increase if we replace $(t x_1 + 1)^{n-2}$ by an appropriate exponential function $ e^{c x_1 + d}$. To this end, we claim that there exist $c, d \in \mathbb{R}$ such that
\begin{equation}\label{1stcond}
\int_0^1 e^{c x_1 + d} dx_1 = \int_0^1 (t x_1 + 1)^{n-2} dx_1 
\end{equation}
and
\begin{equation}\label{2ndcond}
\int_0^1  x_1  e^{c x_1 + d}\, dx_1 = \int_0^1 (t x_1 + 1)^{n-2} x_1 dx_1.
\end{equation}
The integrals in the left hand sides of \eqref{1stcond} and \eqref{2ndcond} are equal to
\begin{align*}
\int_0^1 e^{c x_1 + d} dx_1 = e^d \, \frac{e^c - 1}{c} \quad \textup{and} \quad \int_0^1 x_1 e^{c x_1 + d} dx_1 = e^d \, \frac{ 1 + (c - 1) e^c}{c^2} .
\end{align*}
The value  $c=0$ is included above as a limiting case.
Dividing  \eqref{2ndcond} by \eqref{1stcond}, we get
\begin{align}\label{eq1}
- \frac1{c} + \frac{e^c}{e^c - 1} = \frac{\int_0^1 (t x_1 + 1)^{n-2} x_1 dx_1}{\int_0^1 (t x_1 + 1)^{n-2} dx_1}.
\end{align} 
Observe that the ratio on the right is strictly between $0$ and $1$. 

Consider the function $\phi(c) = - \frac1{c} + \frac{e^c}{e^c - 1}$. We  can assume that $\phi$ is continuous at $c=0$ by setting $\phi(0)=\lim\limits_{c \rightarrow 0} \phi(c) = \frac12.$ Note that  $\lim\limits_{c \rightarrow -\infty} \phi(c) = 0$ and $\lim\limits_{c \rightarrow \infty} \phi(c) = 1$.  Therefore, by continuity of $\phi$,   there exists $c \in \mathbb{R}$ such that \eqref{eq1} holds. Next, choosing    $d \in \mathbb{R}$ so that
$$
e^d = \frac{c}{e^c - 1} \int_0^1 (t x_1 + 1)^{n-2} dx_1
$$
(note that $\frac{e^c - 1}c >0$), guarantees that \eqref{1stcond} and \eqref{2ndcond} are  satisfied.

Now we will  show that 
$$
\int_0^1 (t x_1 + 1)^{n-2} x_1^2 dx_1 \leq \int_0^1 x_1^2 e^{ c x_1 + d} dx_1.
$$
Equation \eqref{1stcond}, which   can be written as follows:
\begin{equation*} 
\int_0^1\left( e^{c x_1 + d} - (t x_1 + 1)^{n-2}\right) dx_1 =0,
\end{equation*} implies  that there is a root $x_1 =\alpha_1 \in (0,1)$ of the function  $\psi(x_1)=e^{c x_1 + d} - (t x_1 + 1)^{n-2}$. Note that such a root cannot be unique. If it were unique, then the function $\psi(x_1)\cdot (x_1-\alpha_1)$ would not change its sign in the interval $(0,1)$. However, by \eqref{1stcond} and  \eqref{2ndcond} we have  $\int_0^1 \psi(x_1) (x_1-\alpha_1)\, dx_1=0$, implying $\psi (x_1)=0$ for all $x_1\in [0,1]$, which is impossible.  On the other hand, the function $\psi$ cannot have more than two roots, since  the exponential function $e^{\frac{c}{n-2} x_1 + \frac{d}{n-2}} $ is convex and its graph cannot intersect the graph of the linear function $ t x_1 + 1$ more than twice. Thus $\psi$ has exactly two roots $\alpha_1$ and $\alpha_2$ in the interval $(0,1)$. We will assume $\alpha_1<\alpha_2$. 
By the convexity of $e^{\frac{c}{n-2} x_1 + \frac{d}{n-2}}$, we have   $e^{\frac{c}{n-2} x_1 + \frac{d}{n-2}} \geq t x_1 + 1$ when $x_1\in [0,\alpha_1]\cup[\alpha_2,1]$, and $e^{\frac{c}{n-2} x_1 + \frac{d}{n-2}} \leq t x_1 + 1$ when $ x_1\in [\alpha_1 , \alpha_2]$.
Using this observation, as well as equalities  \eqref{1stcond} and \eqref{2ndcond}, we get
\begin{align*}
&\int_0^1  \left(e^{c x_1 + d} - (t x_1 + 1)^{n-2}  \right) x_1^2 dx_1 \\
& =\int_0^1  \left(e^{(\frac{c}{n-2} x_1 + \frac{d}{n-2}) (n-2)} - (t x_1 + 1)^{n-2}  \right) (x_1 - \alpha_1) (x_1 - \alpha_2)dx_1 \geq 0.
\end{align*}
 Therefore,
 $$D_n(t) \leq \frac{ \int_0^1 x_1^2 e^{c x_1} dx_1}{\int_0^1 x_1 e^{c x_1} dx_1} - \frac{\int_0^1 x_1 e^{c x_1} dx_1}{\int_0^1  e^{c x_1} dx_1}
 =\frac{c^2 e^c-2ce^c+2e^c-2}{c(1-e^c+ce^c)}-\frac{1-e^c+ce^c}{c(e^c-1)}.$$
Denote $$D= \max_{c\in \mathbb R} \left( \frac{c^2 e^c-2ce^c+2e^c-2}{c(1-e^c+ce^c)}-\frac{1-e^c+ce^c}{c(e^c-1)}\right).$$
The latter expression achieves its maximum at $c_0 \approx - 2.35332$, which can be found numerically. Thus  $D \approx 0.201619.$ 

We have shown that for all $n$ we have 
$$D_n = \sup_{t\ge -1} D_n(t) \le D.$$
To show that $\lim_{n\to \infty} D_n =D$,  consider the sequence $t_n = c_0/n$. Then

\begin{align*}&\lim_{n\to \infty} D_n \ge \lim_{n\to \infty} D_n (t_n)\\
 & =\lim_{n\to \infty}\left( \frac{-2+(t_n+1)^{n-1} (2+(n-1)t_n(nt_n-2))}{(n+1)t_n(1+(t_n+1)^{n-1}((n-1)t_n-1))} -  \frac{1+(t_n+1)^{n-1}((n-1)t_n-1)}{n t_n((t_n+1)^{n-1}-1)}\right)\\
  & =  \frac{-2+e^{c_0} (2+ c_0(c_0-2))}{c_0(1+e^{c_0}(c_0-1))} -  \frac{1+e^{c_0}(c_0-1)}{c_0(e^{c_0}-1)}=D.
  \end{align*}

{\it Step 7 - Classification of all maximizers.}
We will now characterize all possible maximizers.
As we have shown above, for maximizers we have $G(x_1) = \frac{x_1}{n-2}$.
We will now see that this implies  $g(x_1,r)=x_1$. Since $G$ is a linear function,  (\ref{G_concave}) turns into equality, which  implies that $g\left(\tau,  (a \tau +b)z  \right)$ is a linear function of $\tau$. Thus, $g\left(\tau,  (a \tau +b)z  \right)=h_1(z) \tau +h_2(z)$ for some functions $h_1$ and $h_2$ that depend only on $z$.
Observe that $g\ge 0$ implies $g\left(0,  b z  \right)= h_2(z)\ge 0$. Furthermore, since $$
G(\tau) = \int_0^1 z^{n-3} g(\tau,(a   \tau+b)z) \, dz = \frac{\tau}{n-2},
$$
we have 
 $$
0=G(0) = \int_0^1 z^{n-3}h_2(z) \, dz,
$$ and thus $h_2$ is identically equal to zero, i.e.,
$g\left(\tau,  (a \tau +b)z  \right)=h_1(z) \tau.$  Making a change of variables,
we get
$$g(x,y)= h_1\left(\frac{y}{ax+b}\right) x.$$ Since $g(x,y)$ is concave in $y$ and $g(x,y)\le g(x,0)$ for all $(x,y)$ in the domain of $g$, we obtain that $h_1(z)$ is a concave function of $z$ and $h_1(z)\le h_1(0)$ for all  $z\in[0,1]$.

Let $\lambda\in [0,1]$. Using the concavity of $g(x,y)$, we get
$$\lambda g(0,b)+(1-\lambda) g(1,0)\le g(1-\lambda, \lambda b).$$
This implies
$$(1-\lambda)h_1(0)\le (1-\lambda)h_1\left(\frac{\lambda b}{a(1-\lambda)+b}\right),$$
and therefore
$$ h_1(0)=  h_1\left(\frac{\lambda b}{a(1-\lambda)+b}\right),$$
for all $\lambda\in [0,1]$. The latter is equivalent to $h_1(0)=h_1 (z)$ for all $z\in [0,1]$, meaning that $g(x_1,r)=h\cdot  x_1$, for some constant $h$. In our setting this constant is equal to 1, but  this is not important since we can always rescale.

Thus, we have shown that bodies of the form $$K=\{x\in\mathbb R^n:  0\le x_1\le 1,\, \sqrt{x_2^2+\cdots+x_{n-1}^2}\le ax_1+b,\, 0\le x_n\le   x_1\} $$ are the only maximizers in the class $\mathcal C$, with $a$ and $b$ such that $t=a/b$  maximizes the rational function
$D_n(t)$ introduced above. Below we will refer to this $t$ as $t_{max}$. We will now classify all maximizers. Recall that the class $\mathcal C$ was obtained by performing Blaschke shaking and the Schwarz-Steiner symmetrization. We will now undo these operations. We will start by analyzing the projection. Assume $K$ was obtained from a body $  K'\subset \mathbb R^n$ using the Schwarz-Steiner symmetrization with respect to the coordinate plane $x_1x_n$. 
Then, for each $\lambda\in[0,1]$,  by the Brunn-Minkowski inequality we have
\begin{align*}(\kappa_{n-2})^{\frac{1}{n-2}}(a\lambda+b)& =|P_H   K\cap \{x_1=\lambda\}|^{\frac{1}{n-2}}=|P_H  K'\cap \{x_1=\lambda\}|^{\frac{1}{n-2}}\\
&\ge (1-\lambda)|P_H   K'\cap \{x_1=0\}|^{\frac{1}{n-2}}+\lambda|P_H     K'\cap \{x_1=1\}|^{\frac{1}{n-2}}\\
&= (1-\lambda)|P_H    K\cap \{x_1=0\}|^{\frac{1}{n-2}}+\lambda|P_H     K\cap \{x_1=1\}|^{\frac{1}{n-2}}\\
&=(1-\lambda) (\kappa_{n-2})^{\frac{1}{n-2}}b+\lambda(\kappa_{n-2})^{\frac{1}{n-2}}(a+b) = (\kappa_{n-2})^{\frac{1}{n-2}}(a\lambda+b).
\end{align*}
The equality case in the Brunn-Minkowski inequality implies that   $L_0=P_H   K'\cap \{x_1=0\}$ and $L_1=P_H   K'\cap \{x_1=1\}$ are homothetic. Moreover, $P_H  K'$ is the convex hull of $L_0$ and $L_1$.

Next observe that $  K'\cap \{x_1=\lambda\}$ is a cylinder for every $0<\lambda\le 1$. Indeed, for every $0\le s\le \lambda$ we have  
\begin{multline*} |  K'\cap \{x_1=\lambda\}\cap \{x_n=s\}|=| K\cap \{x_1=\lambda\}\cap \{x_n=s\}|\\ =|  K\cap \{x_1=\lambda\}\cap \{x_n=0\}|
=|   K'\cap \{x_1=\lambda\}\cap \{x_n=0\}|.
\end{multline*}
Since $  K'\cap \{x_1=\lambda\}\cap \{x_n=s\}$ and $ K'\cap \{x_1=\lambda\}\cap \{x_n=0\}$ have the same volume, they must be translates of each other by a vector parallel to $e_n$, since otherwise  $P_H (  K'\cap \{x_1=\lambda\})$ would have larger volume than $  K'\cap \{x_1=\lambda\}\cap \{x_n=0\}$.

We conclude that $ K'$ is the convex hull of three $(n-2)$-dimensional bodies $L_0$, $L_1$, and $L_2$, such that $L_0\subset \{x_1=0\}\cap \{x_n=0\}$, $L_1\subset \{x_1=1\}\cap \{x_n=0\}$, $L_2\subset \{x_1=1\}\cap \{x_n=1\}$,  $P_H L_2=L_1$, and $L_0$ is homothetic to $L_1$ (and $L_2$), with $t_{max}$ being the coefficient of homothety. 

Note that for the resulting body $K'$ we need to impose an additional assumption that the segment connecting  $ P_Hc(K')$ and $c(P_HK')$ is  parallel to $e_1$. We will now study this condition. 
Observe that the centroids of all sections $\left(P_H   K'\right)\cap \{x_1=t\}$ lie on the line $\gamma$ connecting the centroids of $L_0$ and $L_1$, and therefore the centroid of $P_H   K'$ also lies on the same line $\gamma$. By the same reasoning, $P_Hc(K')$ belongs to $\gamma$. These observations show that the segment connecting  $ P_Hc(K')$ and $c(P_HK')$ is parallel to $e_1$ if and only if $\gamma$ is parallel to $e_1$.

Thus, the body $K'$ can now be described as the convex hull of $L_0$, $L_1$, and $L_2$, where $L_1$ is obtained from $L_0$ by dilating $L_0$ with a coefficient $t_{max}$ with respect to its centroid and then translating it by $e_1$. $L_2$ is a translate of $L_1$ by   $e_n$.

We will now undo Blaschke shaking.  
Assume that $K'$ is obtained from  $K''$ by applying Blaschke shaking. 
First note that $K'$ and $K''$ have the same projection onto $H$. We can describe the body $K''$ as follows:
$$K'' = \{x\in \mathbb R^n: (x_1,\ldots,x_{n-1})\in P_H K', g_-(x_1,\ldots,x_{n-1}) \le x_n \le g_+(x_1,\ldots,x_{n-1})\}$$
for a convex function $g_-$ and a concave function $g_+$.

Since the length of $\ell \cap K'$ is the same as the length of  $\ell \cap K''$ for every line $\ell$ parallel to $e_n$, we get that 
$$g_+(x_1,\ldots,x_{n-1}) - g_-(x_1,\ldots,x_{n-1})=x_1,$$
which implies that both $g_-$ and $g_+$ are linear functions that differ only in the $x_1$-coordinate. The affine hyperplanes    $x_n=g_-(x_1,\ldots,x_{n-1})$ and $x_n=g_+(x_1,\ldots,x_{n-1})$ intersect $\{x_1=0\}$ in the same $(n-2)$-dimensional affine subspace, which we will denote by $U$. The intersections of the above affine hyperplanes with $\{x_1=1\}$ are  obtained from $U$ by translations. 

Thus $K''$ is the convex hull of three $(n-2)$-dimensional  convex bodies $L_0$, $L_1$, $L_2$ that can be described as follows. Take any $(n-2)$-dimensional subspace $U$ of $\{x_1=0\}$ transversal to $e_n$. Take any convex body $L_0$ in $U$. Let $tL_0$ be the dilation of $L_0$ with respect to its centroid by a factor of $t=t_{max}$.
Let $\mu\ne \nu$ be two real numbers. Define $L_1= tL_0+e_1+\mu e_n$ and $L_2= tL_0+e_1+\nu e_n$.
Dilating  and translating the resulting body $K''$ will give all maximizers of (\ref{MainIneq}) when $H=\{x_n=0\}$ and $u=e_1$. For all other $H$ and $u$ the maximizers obtained by rotations.

 \end{proof}

\section{Projections onto $k$-dimensional subspaces}

It is natural to ask what happens in the case of projections onto lower-dimensional subspaces. Namely, we are interested in
the smallest constant $D_{n,k}$ such that for any convex body $K$ in $\mathbb R^n$ and every $k$-dimensional subspace $H$ we have $$ |P_Hc(K)-c(P_HK)|\le D_{n,k} w_K (u),$$ where $u$ is the unit vector parallel to the segment connecting $ P_Hc(K)$ and $c(P_HK)$.
 
 To solve this problem, one can use methods similar to those above. Namely, assuming that $H=\mbox{span}\{e_1,\ldots,e_k\}$, $u=e_1$, and  that   $\{x_1=0\}$ and $\{x_1=1\}$ are supporting hyperplanes to the body $K$, one can show that the maximum of $|P_Hc(K)-c(P_HK)|$ is achieved for bodies of the form
 $$\{x\in\mathbb R^n:0\le x_1\le 1,\, \sqrt{x_2^2+\cdots+x_k^2}\le ax_1+b,\, \sqrt{x_{k+1}^2+\cdots+x_n^2} \le x_1\}.$$
 From this we can get
  \begin{align*}
D_{n,k}=\sup_{t\ge -1}\left( \frac{\int_0^1 (t u + 1)^{k-1} u^{n-k+1}\, du}{\int_0^1 (t u + 1)^{k-1} u^{n-k}\,  du} - \frac{\int_0^1 (t u + 1)^{k-1} u \,  du}{\int_0^1 (t u + 1)^{k-1}  \,  du}\right).
 \end{align*}
 Unfortunately, as in the previous section, it is impossible to compute $D_{n,k}$ precisely for all $n$ and $k$. All we can do is to show that for a fixed $k$ the sequence $D_{n,k}$ is increasing with respect to $n$ and $$\lim_{n\to\infty} D_{n,k} = 1 - \frac{1}{k+1}.$$
 Therefore, $$ |P_Hc(K)-c(P_HK)|\le \left(1-\frac{1}{k+1}\right) w_K (u).$$
 We will however omit the details, since this bound can be obtained much easier. Indeed, clearly in the class of bodies defined above we have
 $\langle P_Hc(K),e_1\rangle \le 1$ and $\langle c(P_HK),e_1\rangle \ge \frac{1}{k+1}$. The latter is the well-known fact that among all convex bodies with centroid at the origin the minimum of  $ \max_{x\in K} \{\langle x,e_1\rangle\}/ \max_{x\in K} \{\langle x,-e_1\rangle\}$ is achieved for cones 
 whose base is orthogonal to the $x_1$-axis; cf. \cite[p.~58]{BF}.

\bigskip

{\bf Acknowledgment.} The authors are grateful to Fedor Nazarov  for his invaluable help.

\end{document}